%% file: arx3.tex
\documentclass[11pt]{amsart}
\usepackage{url,graphicx}
\usepackage[dvips]{epsfig}
\usepackage{latexsym,color}
\usepackage{amscd,amssymb,verbatim}
\newcommand{\be}{\begin{enumerate}}
\newcommand{\ee}{\end{enumerate}}

\newcommand{\da}{\Delta}
\newcommand{\thom}{\text{\tt Hom}\,}
\newcommand{\ra}{\rightarrow}
\newcommand{\nin}{\noindent}
\newcommand{\pr}{\noindent{\bf Proof. }}
\newcommand{\sm}{\setminus}
\newcommand{\dne}{\searrow_{\text{\tt NE}}}
\newcommand{\une}{\nearrow_{\text{\tt NE}}}
\newcommand{\lk}{\text{\rm lk}}
\newcommand{\ene}{\simeq_{\text{\tt NE}}}
\newcommand{\noe}{{\tt NE}}
\newcommand{\fix}{\text{\rm Fix}\,}

\newtheorem{thm}{Theorem}[section]
\newtheorem{df}[thm]{Definition}

\newtheorem{crl}[thm]{Corollary}
\newtheorem{prop}[thm]{Proposition}
\newtheorem{conj}[thm]{Conjecture}

\numberwithin{equation}{section}
\numberwithin{figure}{section}

\begin{document}

\title
{Collapsing along monotone poset maps}

\author{Dmitry N. Kozlov}
\address{Institute of Theoretical Computer Science, Eidgen\"ossische Technische
Hochschule, Z\"urich, Switzerland}
\email{dkozlov@inf.ethz.ch}
\thanks {Research supported by Swiss National Science Foundation Grant PP002-102738/1}

\date\today

\begin{abstract}
We introduce the notion of nonevasive reduction, and show that for any monotone poset map $\varphi:P\ra P$, the simplicial complex $\da(P)$ \noe-reduces to $\da(Q)$, for any $Q\supseteq\fix\varphi$. 

As a~corollary, we prove that for any order-preserving map $\varphi:P\ra P$ satisfying $\varphi(x)\geq x$, for any $x\in P$, the simplicial complex $\da(P)$ collapses to $\da(\varphi(P))$. We also obtain a~generalization of Crapo's closure theorem.  
\end{abstract}

\maketitle

\section{Order complexes, collapsing and \noe-reduction.}

For a poset $P$ we let $\Delta(P)$ denote its {\it nerve}: the
simplicial complex whose simplices are all chains of $P$. For a
simplicial complex $X$ we let $V(X)$ denote the set of its vertices.

An {\it elementary collapse} in a~simplicial complex $X$ is a~removal
of two open simplices $\sigma$ and $\tau$ from $X$, such that
$\dim\sigma=\dim\tau+1$, and $\sigma$ is the only simplex of $X$,
different from $\tau$ itself, which contains the simplex $\tau$ in its
closure.  

When $Y$ is a~subcomplex of $X$, we say that $X$ {\it collapses onto}
$Y$ if there exists a~sequence of elementary collapses leading from
$X$ to $Y$; in this case we write $X\searrow Y$ (or, equivalently,
$Y\nearrow X$).

\begin{df}\label{df:nonev} $\,$

\nin (1) A finite nonempty simplicial complex $X$ is called {\bf nonevasive} if either $X$ is a~point, or, inductively, there exists a vertex $v$ of $X$, such that both $X\setminus\{v\}$ and $\lk_X v$ are nonevasive.

\nin (2) For two nonempty simplicial complexes $X$ and $Y$ we write $X\dne Y$ (or, equivalently, $Y\une X$), if there exists a~sequence $X=A_1\supset A_2\supset\dots\supset A_t=Y$, such that for all $i\in\{1,\dots,t-1\}$ we can write $V(A_i)=V(A_{i+1})\cup\{x_i\}$, so that $lk_{A_i}x_i$ is nonevasive.
\end{df}

We recall, that the notion of nonevasive simplicial complexes was
introduced in \cite{KSS}, and was initially motivated by the
complexity-theoretic considerations. For further connections to
topology and more facts on nonevasiveness we refer to \cite{Ku,We}.
Recently an~interesting connection has been established between
Discrete Morse theory and evasiveness, the standard references are
\cite{F1,F2}.

Several classes of simplicial complexes are known to be
nonevasive. Perhaps the simplest example is provided by the fact that
all cones are nonevasive. A~more complicated family of nonevasive
simplicial complexes is obtained by taking the order complexes of the
noncomplemented lattices, see~\cite{Ko2}.

In the situation described in Definition \ref{df:nonev}(2), we say
that the simplicial complex $X$ \noe-reduces to its subcomplex
$Y$. The following facts about \noe-reduction are useful for our
arguments.

\vspace{5pt}

\nin {\bf Fact 1.}
{\it If $X_1$ and $X_2$ are simplicial complexes, such that $X_1\dne X_2$, and $Y$ is an arbitrary simplicial complex, then $X_1*Y\dne X_2*Y$.} 

\vspace{5pt}

\nin Here the symbol $*$ denotes the simplicial join of two simplicial 
complexes, see \cite{Munk}. That $X_1*Y\dne X_2*Y$ follows by
induction from the facts that if $v$ is any vertex of a~simplicial
complex $X$, then we have $\lk_{X*Y}v=(\lk_{X}v)*Y$ and
$(X*Y)\sm\{v\}=(X\sm\{v\})*Y$, for any $v\in V(X)$. 

To \noe-reduce $X_1*Y$ to $X_2*Y$, simply take the sequence of
vertices $x_1,\dots,x_t\in V(X_1)$ which \noe-reduces $X_1$ to
$X_2$. We have $\lk_{X_1*Y}(x_1)=(\lk_{X_1}(x_1))*Y$. In turn, the
simplicial complex $(\lk_{X_1}(x_1))*Y$ is nonevasive: this is seen by
induction on the number of vertices of the first factor, with the base
given by the fact that all cones are nonevasive. Removing $x_1$ from
$X_1*Y$ yields the simplicial complex $(X_1\sm\{x_1\})*Y$, hence,
continuing in this way, we will \noe-reduce $X_1*Y$ all the way to
$X_2*Y$.

\vspace{5pt}

\nin {\bf Fact 2.}
{\it The reduction $X\dne Y$ implies $X\searrow Y$, which in turn
implies that $Y$ is a~strong deformation retract of~$X$.}

\section{Monotone poset maps.}

Next we define a~class of maps which are particularly suitable for our purposes.

\begin{df} \label{df:cpa}
Let $P$ be a poset. An order-preserving map $\varphi:P\ra P$ is called
a~{\bf monotone map}, if for any $x\in P$ either $x\geq\varphi(x)$ or
$x\leq\varphi(x)$.

If $x\geq\varphi(x)$ for all $x\in P$, then we call $\varphi$ a~{\bf
decreasing map}, analogously, if $x\leq\varphi(x)$ for all $x\in P$,
then we call $\varphi$ an~{\bf increasing map}.
\end{df}

We remark here the fact that while a composition of two decreasing
maps is again a~decreasing map, and, in the same way, a~composition of
two increasing maps is again an~increasing map, the composition of two
monotone maps is not necessarily a~monotone map. 

\vskip5pt

\nin {\bf Example.} {\it Let $P$ be the lattice of all subsets of $\{1,2\}$, 
and define $\varphi(S)=S\cup\{2\}$, and $\gamma(T)=T\sm\{1\}$, for all
$S,T\subseteq\{1,2\}$. The composition $T\circ S$ maps all the subsets
to $\{2\}$, in particular it is not a~monotone map.}

\vskip5pt

On the other hand, any power of a~monotone map is again monotone. Indeed, let $\varphi:P\ra P$ be monotone, let $x\in P$, and say $x\leq\varphi(x)$. Since $\varphi$ is order-preserving we conclude that $\varphi(x)\leq\varphi^2(x)$, $\varphi^2(x)\leq\varphi^3(x)$, etc. Hence $x\leq\varphi^N(x)$ for arbitrary~$N$.

The following proposition shows that monotone maps have a canonical decomposition in terms of increasing and descreasing maps.

\begin{prop} \label{pr:mondec}
Let $P$ be a poset, and let $\varphi:P\ra P$ be a~monotone map. There exist unique maps $\alpha,\beta:P\ra P$, such that 
\begin{itemize}
	\item $\varphi=\alpha\circ\beta$;
	\item $\alpha$ is an increasing map, whereas $\beta$ is a~decreasing map;
	\item $\fix\alpha\cup\fix\beta=P$.
\end{itemize}
\end{prop}

\pr Set 
\begin{equation}\label{aleq}
\alpha(x):=\begin{cases}
\varphi(x),& \text{ if } \varphi(x)>x;\\
x,& \text{ otherwise,}
\end{cases}
\end{equation}
and
\begin{equation}\label{beteq}
\beta(x):=\begin{cases}
\varphi(x),& \text{ if } \varphi(x)<x;\\
x,& \text{ otherwise.}
\end{cases}
\end{equation}
Clearly, $\fix\alpha\cup\fix\beta=P$. Let us see that
$\varphi(x)=\alpha(\beta(x))$, for any $x\in P$. This is obvious if
$\varphi(x)\geq x$, since then $\alpha(\beta(x))=\alpha(x)=\varphi(x)$ by
\eqref{aleq} and \eqref{beteq} respectively. Assume $\varphi(x)<x$, then
$\beta(x)=\varphi(x)$, hence
$\alpha(\beta(x))=\alpha(\varphi(x))$. Since $\varphi$ is
order-preserving, $\varphi(x)<x$ implies
$\varphi(\varphi(x))\leq\varphi(x)$. Thus \eqref{aleq} gives
$\varphi(x)\in\fix\alpha$, and we conclude that
$\alpha(\beta(x))=\alpha(\varphi(x))=\varphi(x)$.

To see that $\alpha$ is an~increasing map, we just need to see that it
is order-preserving. Since $\alpha$ either fixes an element or maps it
to a~larger one, the only situation which needs to be considered is
when $x,y\in P$, $x<y$, and $\alpha(x)=\varphi(x)$. However, under
these conditions we have
$\alpha(y)\geq\varphi(y)\geq\varphi(x)=\alpha(x)$, and so $\alpha$ is
order-preserving. That $\beta$ is a~decreasing map can be seen
analogously. Finally, the uniqueness follows from the fact that each
$x\in P$ must be fixed by either $\alpha$ or~$\beta$, and the value
$\varphi(x)$ determines which one will fix~$x$.
\qed

\section{The main theorem and implications.}

Prior to this work, it has been known that a monotone map $\varphi:P\ra P$ induces a~homotopy equivalence between $\da(P)$ and $\da(\varphi(P))$, see \cite[Corollary 10.17]{Bj}. It was also proved in \cite{Bj} that if the map $\varphi$ satisfies the additional condition $\varphi^2=\varphi$, then $\varphi$ induces a~strong deformation retraction from $\da(P)$ to $\da(\varphi(P))$. 

The latter result was strengthened in \cite[Theorem 2.1]{Ko1}, where it was shown that, whenever $\varphi$ is an~ascending (or descending) closure operator, $\da(P)$ collapses onto $\da(\varphi(P))$. There this fact was used to analyze the effect of the folding operation on the corresponding $\thom$ complexes, see also \cite{BK03a,BK03b,BK03c,IAS}.

The next theorem strengthens and generalizes the results from \cite{Bj} and \cite{Ko1}.

\begin{thm} \label{thm1}
Let $P$ be a~poset, and let $\varphi:P\ra P$ be a~monotone map.
Assume $P\supseteq Q\supseteq\fix\varphi$, $P\sm Q$ is finite, and,
for every $x\in P\sm Q$, $P_{<x}\cup P_{>x}$ is finite, then
$\Delta(P)\dne\Delta(Q)$, in particular, $\Delta(P)$ collapses onto
$\Delta(Q)$.
\end{thm}

\nin {\bf Remarks.}

\vskip5pt

\nin 1) Note that  when $P$ is finite, the conditions of the Theorem \ref{thm1} simply reduce to: $P\supseteq Q\supseteq\fix\varphi$.

\vskip5pt

\nin 2) Under conditions of Theorem \ref{thm1}, the simplicial complex $\Delta(P)$ collapses onto the simplicial complex $\Delta(Q)$, in particular, the complexes $\da(P)$ and $\da(Q)$ have the same simple homotopy type, see \cite{Co}.

\vskip5pt

\nin 3) Under conditions of Theorem \ref{thm1}, the topological space $\da(Q)$ is a strong deformation retract of the topological space $\da(P)$.

\vskip5pt

\nin 4) Any poset $Q$ satisfying $P\supseteq Q\supseteq\varphi(P)$ will also
satisfy $P\supseteq Q\supseteq\fix\varphi$, hence Theorem~\ref{thm1}
will apply. In particular, for finite $P$, we have the following
corollary:

\begin{crl}
For any poset $P$, and for any monotone map $\varphi:P\ra P$,
satisfying conditions of Theorem~\ref{thm1}, we have
$\da(P)\dne\da(\varphi(P))$.
\end{crl}

It is easy to prove Theorem \ref{thm1}, once the following auxiliary
result is established.

\begin{prop} \label{prop3.1a}
Let $P$ be a~poset, and let $\varphi:P\ra P$ be a~monotone map.
Assume $x\in P$, such that $\varphi(x)\neq x$, and $P_{<x}\cup P_{>x}$
is finite, then $\da(P_{<x})*\da(P_{>x})$ is nonevasive. More precisely, if $\varphi(x)<x$, then $\da(P_{<x})$ is nonevasive, and if $\varphi(x)>x$, then $\da(P_{>x})$ is nonevasive.
\end{prop}

\pr Since the expression $\da(P_{<x})*\da(P_{>x})$ is symmetric with respect to inverting the partial order of $P$, it is enough, without loss of generality, to only consider the case $\varphi(x)< x$. Let us show that in this case $\da(P_{<x})$ is nonevasive. We proceed by induction on $|P_{<x}|$. If $|P_{<x}|=1$, then the statement is clear, so assume $|P_{<x}|\geq 2$.

Let $\psi:P_{<x}\ra P_{<x}$ denote the restriction of $\varphi$. It is easy to see that $\psi$ is a~monotone map of $P_{<x}$. To verify that $\Delta(P_{<x})\dne\Delta(P_{\leq\varphi(x)})$, order the elements of $P_{<x}\sm P_{\leq\varphi(x)}$ following an~arbitrary linear extension in the decreasing order, say $P_{<x}\sm P_{\leq\varphi(x)}=\{a_1,\dots,a_t\}$, and $a_i\not< a_j$, for $i<j$. By the choice of the order of $a_i$'s, we have $P_{<a_i}=P^i_{<a_i}$, where $P^i=P\sm\{a_1,\dots,a_{i-1}\}$. Therefore, by the induction assumption, $\da(P_{<a_i})$ is nonevasive for all $1\leq i\leq t$, and we have \[\da(P_{<x})=\da(P^1_{<x})\dne\da(P^2_{<x})\dne\dots\dne\da(P^{t+1}_{<x})=\Delta(P_{\leq\varphi(x)}).\]

On the other hand, $\Delta(P_{\leq\varphi(x)})$ is a cone, hence it is nonevasive, and therefore $\Delta(P_{<x})$ is nonevasive as well. It follows that $\da(P_{<x})*\da(P_{>x})$ is nonevasive.
\qed

\nin {\bf Proof of the Theorem \ref{thm1}.} The proof is by induction on $|P\sm Q|$. The statement is trivial when $|P\sm Q|=0$, so assume $|P\sm Q|\geq 1$.

To start with, we replace the monotone map $\varphi$ with a~monotone
map $\gamma$ satisfying $\gamma(P)\subseteq Q$ and
$\fix\gamma=\fix\varphi$. To achieve that objective we can set
$\gamma:=\varphi^N$, where $N=|P\sm Q|$. With this choice of $\gamma$,
the inclusion $\gamma(P)\subseteq Q$ follows from the assumption that
$\fix\varphi\subseteq Q$, since $\fix\gamma=\gamma(P)$.

Take arbitrary $x\in P\sm Q$. Since $x\notin\gamma(P)$, we have $x\neq\gamma(x)$, hence by Proposition~\ref{prop3.1a} we know that $\lk_{\da(P)}x=\da(P_{<x})*\da(P_{>x})$ is nonevasive. This means $\da(P)\dne\da(P\sm\{x\})$. 

Let $\psi:P\sm\{x\}\ra P\sm\{x\}$ be the restriction of $\gamma$. Clearly, $\psi$ is a~monotone map, and $\fix\psi=\fix\gamma$. This implies that $\fix\psi\subseteq Q$, hence, by the induction hypothesis $\da(P\sm\{x\})\dne\da(Q)$. Summarizing, we conclude that $\da(P)\dne\da(Q)$.
\qed

\vskip5pt

On the enumerative side, we obtain the following generalization of the Crapo's Closure Theorem from 1968, see \cite[Theorem 1]{Cr}.

\begin{crl} \label{crl:cr}
Let $P$ be a finite poset with $\hat 0$ and $\hat 1$, and let $\varphi:P\ra P$ be an~increasing map. Assume $P\supseteq Q\supseteq\fix\varphi$, and $Q\cap\varphi^{-\infty}(\hat 1)=\{\hat 1\}$. Then, we have
\[
\sum_{\varphi^\infty(z)=\hat 1}\mu_P(\hat 0,z)=
\begin{cases}
\mu_Q(\hat 0,\hat 1),&\text{ if } \hat 0\in\fix\varphi;\\
0,&\text{ otherwise}.
\end{cases}
\]
Here $\varphi^\infty$ is the stabilization of $\varphi$, say $\varphi^\infty:=\varphi^{|P|}$, so $\varphi^{-\infty}(\hat 1)$ denotes the set of all elements of $P$ which map to $\hat 1$ after a~sufficiently high iteration of~$\varphi$.
\end{crl}

Before we give the proof, recall the following convention: whenever $P$ is a~poset with $\hat 0$ and $\hat 1$, we let $\bar P$ denote $P\sm\{\hat 0,\hat 1\}$.

\vskip5pt

\nin {\bf Proof of the Corollary \ref{crl:cr}.} Assume first that $\hat 0\in\fix\varphi$, hence $\hat 0,\hat 1\in Q$, and, since $\varphi$ is increasing, $\varphi^{-\infty}(\hat 0)=\{\hat 0\}$. Set $R:=(P\sm\varphi^{-\infty}(\hat 1))\cup\{\hat 0,\hat 1\}$, i.e., $R$ is the set of all elements of $\bar P$ which shall not map to $\hat 1$, no matter how high iteration of $\varphi$ we take, with the maximal and the minimal elements attached. Let $\psi:\bar R\ra\bar R$ be the restriction of $\varphi$. Clearly, $\psi$ is a~monotone map, and $\fix\psi= \fix\varphi\sm\{\hat 0,\hat 1\}$. 

Since $\bar R\supseteq\bar Q\supseteq\fix\psi$, we conclude that $\da(\bar R)$ collapses onto $\da(\bar Q)$; in particular the simplicial complexes $\da(\bar R)$ and $\da(\bar Q)$ have the same Euler characteristic. By Ph.\ Hall Theorem, see \cite{St}, for any poset $P$ with a~maximal and a~minimal element we have $\chi(\da(\bar P))=\mu_P(\hat 0,\hat 1)$, therefore here we conclude that $\mu_Q(\hat 0,\hat 1)=\mu_R(\hat 0,\hat 1)$. 

On the other hand, by definition of the M\"obius function, we have the equality $\sum_{z\in P}\mu_P(\hat 0,z)=0$, which can be rewritten as $\sum_{z\in\varphi^{-\infty}(\hat 1)}\mu_P(\hat 0,z)= -\sum_{z\notin\varphi^{-\infty}(\hat 1)}\mu_P(\hat 0,z)$. Similarly, $\mu_R(\hat 0,\hat 1)=-\sum_{x\in R\sm\{\hat 1\}}(\hat 0,x)$. Since the condition $z\notin\varphi^{-\infty}(\hat 1)$ is equivalent to the condition $z\in R\sm\{\hat 1\}$, and $R_{<x}=P_{<x}$, for any $x\in R$, we conclude that $\mu_Q(\hat 0,\hat 1)=\mu_R(\hat 0,\hat 1)=\sum_{z\in\varphi^{-\infty}(\hat 1)}\mu_P(\hat 0,z)$.

Consider now the case $\hat 0\notin\fix\varphi$. If $\varphi^\infty(\hat 0)=\hat 1$, then the statement follows from the definition of the M\"obius function, since then $\varphi^\infty(z)=\hat 1$ for all $z\in P$. Assume $\varphi^\infty(\hat 0)\neq\hat 1$. We can define a new map $\psi$ by changing the value of $\varphi$ in one element:
	\[\psi(x)=\begin{cases}
	\hat 0, & \text{ if } x=\hat 0;\\
	\varphi(x), & \text{ otherwise.}
\end{cases}
\]
Clearly, $\psi$ is a monotone function, $\fix\psi=\fix\varphi\cup\{\hat 0\}$, and $\varphi^{-\infty}(\hat 1)=\psi^{-\infty}(\hat 1)$. Hence, the first part of the proof applies, and we conclude that 
\[\sum_{\psi^\infty(z)=\hat 1}\mu_P(\hat 0,z)=\sum_{\varphi^\infty(z)=\hat 1}\mu_P(\hat 0,z)=\mu_Q(\hat 0,\hat 1),\] 
for any $Q$ such that $P\supseteq Q\supseteq\fix\psi$, and such that $Q\cap\varphi^{-\infty}(\hat 1)=\{\hat 1\}$. Choose $Q=(P_{\geq\varphi(\hat 0)}\sm\varphi^{-\infty}(\hat 1))\cup\{\hat 0,\hat 1\}$. Since $\varphi(\hat 0)\notin\varphi^{-\infty}(\hat 1)$, the poset $Q$ has only one atom $\varphi(\hat 0)$, thus we have $\mu_Q(\hat 0,\hat 1)=0$, and the proof is complete.
 \qed

\section{\noe-reduction and collapses}

The \noe-reduction can be used to define an interesting equivalence relation on the set of all simplicial complexes.

\begin{df} \label{df:noe}
Let $X$ and $Y$ be simplicial complexes. Recursively, we say that $X\ene Y$ if $X\dne Y$, or $Y\une X$, or if there exists a~simplicial complex $Z$, such that $X\ene Z$ and $Y\ene Z$.
\end{df}

\begin{figure}[hbt]
\begin{center}
  \begin{picture}(0,0)%
    \includegraphics{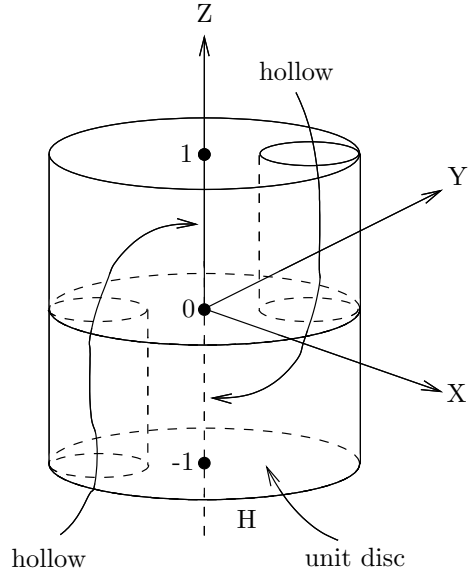}%
  \end{picture}%
  \input{h2r.pstex_t}%
  
\end{center}
\caption{A house with two rooms.}
\label{fig:h2r}
\end{figure}

Clearly, if $X$ is nonevasive, then $X\ene\text{pt}$, but is the
opposite true? The answer to that is "no". To see this, consider the
standard example of a~space which is contractible, but not
collapsible: let $H$ be the so-called {\it house with two rooms}, see
Figure \ref{fig:h2r}.

The space $H$ is not collapsible, hence nonevasive, see \cite{Co} for an argument. On the other hand, we leave it to the reader to see that it is possible to triangulate the filled cylinder $C$ given by the equations $|z|\leq 1$, $x^2+y^2\leq 1$, so that $C\dne H$.

The analogous equivalence relation, where $\dne$, and $\une$, are
replaced by $\searrow$, and $\nearrow$, is called the simple homotopy
equivalence; its equivalence classes are called simple homotopy
types. The celebrated theorem of J.H.C.\ Whitehead states that the
simplicial complexes with the simple homotopy type of a~point are
precisely those, which are contractible, see~\cite{Co}. Therefore, the
class of the simplicial complexes which are $\noe$-equivalent to
a~point relates to nonevasiveness in the same way as contractibility
refers to collapsibility. Clearly, this means that this class should
constitute an interesting object of study.

We conjecture that the $\noe$-equivalence is much finer than the
Whitehead's simple homotopy type. We make two conjectures: a~weak and
a~strong one.

\begin{conj} \label{conj:noesh1}
There exist finite simplicial complexes $X$ and $Y$ having the same
simple homotopy type, such that $X\not\ene Y$.
\end{conj}

\begin{conj} \label{conj:noesh2}
There exists an infinite family of finite simplicial complexes
$\{X_i\}_{i=1}^\infty$, which all have the same simple homotopy type,
such that $X_i\not\ene X_j$, for all $i\neq j$.
\end{conj}

Again, in the simple homotopy setting, the phenomenon of the
Conjectures \ref{conj:noesh1} and \ref{conj:noesh2} is governed by
an~algebraic invariant called the {\it Whitehead torsion}, namely:
a~homotopy equivalence between finite connected CW-complexes is simple
if and only if its Whitehead torsion is trivial,
see~\cite[(22.2)]{Co}. It is enticing to hope for an existence of some
similar invariant in our
\noe-setting.

Finally, let us remark, that whenever we have simplicial complexes
$X\ene Y$, there exists a~simplicial complex $Z$, such that $X\une
Z\dne Y$. Indeed, assume $A\dne B\une C$, for some simplicial
complexes $A,B$, and $C$. Let $S=V(A)\sm V(B)$, and $T=V(C)\sm
V(B)$. Let $D$ be the simplicial complex obtained by attaching to $A$
the vertices from $T$ in the same way as they would be attached to
$B\subseteq A$. Clearly, since the links of the vertices from $S$ did
not change, they can still be removed in the same fashion as before,
and therefore we have $A\une D\dne C$. Repeating this operation
several times, and using the fact that the reductions $\une$ (as well
as $\dne$) compose, we prove the claim.

\vskip10pt

\nin {\bf Acknowledgments.} We would like to thank the Swiss National 
Science Foundation and ETH-Z\"urich for the financial support of this
research. We also thank the referee for several valuable comments.

\end{document}

%% file: h2r.pstex_t
\begin{picture}(0,0)%
\includegraphics{h2r.pstex}%
\end{picture}%
\setlength{\unitlength}{3947sp}%
\begingroup\makeatletter\ifx\SetFigFont\undefined%
\gdef\SetFigFont#1#2#3#4#5{%
  \reset@font\fontsize{#1}{#2pt}%
  \fontfamily{#3}\fontseries{#4}\fontshape{#5}%
  \selectfont}%
\fi\endgroup%
\begin{picture}(2870,3530)(1107,-2811)
\put(2901,254){\makebox(0,0)[b]{\smash{{\SetFigFont{10}{12.0}{\rmdefault}{\mddefault}{\updefault}{\color[rgb]{0,0,0}hollow}%
}}}}
\put(2591,-2561){\makebox(0,0)[b]{\smash{{\SetFigFont{10}{12.0}{\rmdefault}{\mddefault}{\updefault}{\color[rgb]{0,0,0}H}%
}}}}
\put(2211,-266){\makebox(0,0)[b]{\smash{{\SetFigFont{10}{12.0}{\rmdefault}{\mddefault}{\updefault}{\color[rgb]{0,0,0}1}%
}}}}
\put(1346,-2811){\makebox(0,0)[b]{\smash{{\SetFigFont{10}{12.0}{\rmdefault}{\mddefault}{\updefault}{\color[rgb]{0,0,0}hollow}%
}}}}
\put(2326,614){\makebox(0,0)[b]{\smash{{\SetFigFont{10}{12.0}{\rmdefault}{\mddefault}{\updefault}{\color[rgb]{0,0,0}Z}%
}}}}
\put(3271,-2806){\makebox(0,0)[b]{\smash{{\SetFigFont{10}{12.0}{\rmdefault}{\mddefault}{\updefault}{\color[rgb]{0,0,0}unit disc}%
}}}}
\put(3916,-411){\makebox(0,0)[b]{\smash{{\SetFigFont{10}{12.0}{\rmdefault}{\mddefault}{\updefault}{\color[rgb]{0,0,0}Y}%
}}}}
\put(3911,-1771){\makebox(0,0)[b]{\smash{{\SetFigFont{10}{12.0}{\rmdefault}{\mddefault}{\updefault}{\color[rgb]{0,0,0}X}%
}}}}
\put(2226,-1236){\makebox(0,0)[b]{\smash{{\SetFigFont{10}{12.0}{\rmdefault}{\mddefault}{\updefault}{\color[rgb]{0,0,0}0}%
}}}}
\put(2191,-2201){\makebox(0,0)[b]{\smash{{\SetFigFont{10}{12.0}{\rmdefault}{\mddefault}{\updefault}{\color[rgb]{0,0,0}-1}%
}}}}
\end{picture}%